\documentclass[a4paper,12pt]{article}
\usepackage{amsmath,amsfonts,amsthm}
\usepackage{euscript}
\usepackage{mathbbol}
\usepackage{natbib}
\begin{document}
\title{An Algorithmic and a Geometric Characterization of Coarsening
  At Random}

\author{Richard D. Gill\thanks{URL: {\tt www.math.leidenuniv.nl/{\~{}}gill}.
    Also affiliated with EURANDOM, Eindhoven, the Netherlands, {\tt
      www.eurandom.nl}.}  \\ Leiden University  \\ The
  Netherlands \and Peter D.  Gr\"unwald\thanks{URL: {\tt www.grunwald.nl}.
Also affiliated with
    EURANDOM, Eindhoven, the Netherlands, {\tt www.eurandom.nl}. CWI is the National Research Institute for Mathematics and Computer Science in the Netherlands.} \\ CWI, Amsterdam \\ The Netherlands} \date{\today}

\newtheorem{theorem}{Theorem}
\newtheorem{example}{Example}
\newtheorem{question}{Question}
\newtheorem{definition}{Definition}
\newtheorem{remark}{Remark}
\newtheorem{fact}{Fact}
\renewcommand{\vec}{\ensuremath{\mathbf}}
\def\boldpi{{\mbox{\boldmath$\pi$}}}
\def\boldlambda{{\mbox{\boldmath$\lambda$}}}
\maketitle

\begin{abstract}
We show that the class of conditional distributions satisfying the
coarsening at Random (CAR) property  for discrete data has a simple 
and robust algorithmic description based on randomized uniform multicovers:
combinatorial objects generalizing the notion of partition of a set.
However, the complexity of a given CAR mechanism can be large:
the maximal  ``height''  of the needed multicovers 
can be exponential in the number of points in the sample space.  
The results stem from a geometric
interpretation of the set of CAR distributions as a convex polytope
and a characterization of its extreme points. The hierarchy of CAR
models defined in this way could be useful in parsimonious
statistical modelling of CAR mechanisms, though the results 
also raise doubts in applied work as to the meaningfulness
of the CAR assumption in its full generality.
\end{abstract}
\ 
\newline
\\
{\em This paper has been accepted for publication in the Annals of
  Statistics. It will appear later in 2007 or in 2008.}

\section{Introduction}

In statistical practice one is often presented with incomplete, or
more generally, {\em coarse\/} data. To properly model such data, one
needs to take into account the mechanism by which the data are
coarsened. In practice the details of this coarsening mechanism are
often unknown or computationally expensive to model. Therefore, it is
of interest to determine conditions under which this mechanism can be
safely ignored.  The ``coarsening at random'' (CAR) assumption is the
weakest condition giving this guarantee. It was identified by
\cite{HeitjanR91}. More recently, \cite{GrunwaldH03} and \cite{Jaeger05b}
have stressed that the importance of CAR is not restricted to
statistical applications: when updating a probability distribution
based on new information in learning, artificial intelligence,
or other scientific applications, it precisely characterizes when one can
ignore the distinction between the fact that an event has been
\emph{observed}, and the fact that an event has
\emph{happened}, thereby considerably
simplifying the update process.

Thus, both in statistical inference with coarsened data
and for probability updating in learning algorithms, 
it is attractive to be able to make the CAR assumption. 
In order to be able to judge whether or not the assumption is
warranted, it is important to fully understand its meaning.
Here we approach this problem by giving two intimately
related characterizations of the CAR assumption.
First, we show that the set
of all CAR mechanisms for a given finite sample space can be seen as a
convex polytope. Each CAR mechanism is a mixture of CAR mechanisms
which correspond to the vertices of the polytope. Our first main
result, Theorem~\ref{thm:geometric}, characterizes these vertices.
Our second result, which follows easily from the first, 
complements this geometric view with an algorithmic one. We show that a
simple probabilistic algorithm can simulate any possible CAR
mechanism, and only CAR mechanisms. Prompted by \cite{GillLR97}, earlier
authors \citep{GrunwaldH03,Jaeger05b} have also searched for such
constructions, calling them \emph{procedural models} for CAR.  Yet the
procedural models proposed so far are not quite satisfactory, because
in all cases,
\begin{enumerate} 
\item The procedural model depends on parameters which have to be
fine-tuned in order to guarantee the CAR property; or equivalently, 
\item A small
perturbation in the parameters can destroy the CAR property. 
\end{enumerate}
This ``frailty'' or lack of robustness is an indication that such
procedures may not occur naturally. In fact \cite[Theorem
4.17]{Jaeger05b} shows that the only CAR mechanisms which a
robust procedure can generate must be of a special type known as
``coarsening completely at random'', CCAR.

Here we present a natural way to generate 
all CAR mechanisms, and only CAR mechanisms, 
\emph{that does not require fine-tuning of parameters}.
Our algorithm works for arbitrary finite sample spaces. 
It is based on a generalization 
of the notion of a partition of a set which we call a
\emph{uniform multicover}, or just multicover for short.

Superficially, its existence would have to contradict Jaeger's
theorem mentioned above. But of course, a proven 
theorem does not allow any contradictions. The difference
lies in the notion we use of robustness 
and of its negation, frailty.
Our result can be seen as criticism of Jaeger's
notion of robustness, 
even though this does at first-sight seem appealing
and natural. By parameterizing CAR
distributions in a different manner,
we obtain a representation in which CAR 
\emph{can} be generated without parameter tuning. 
In a nut-shell: we consider a discrete uniform distribution
to be a robust and natural object. Jaeger considers it
to be an easily perturbed object.

We emphasize that  the body of Jaeger's work remains highly relevant;
this is just one of a number of important results he has obtained, and
we too come to the conclusion that CAR mechanisms which are
not CCAR will be very rare in practice.
For instance, our final result, Theorem~\ref{thm:fibo}, shows
that, although no fine-tuning is needed, the complexity (defined in
terms of the ``height'' of multicovers) of the CAR mechanisms
generated by our algorithm can grow exponentially in the size of the
sample space.

The paper is organized as follows. In Section~\ref{sec:prelim} we
briefly introduce coarsening at random and other preliminaries.  In
Section~\ref{sec:geometric} we give our geometric interpretation of
CAR distributions (Theorem~\ref{thm:geometric}). In
Section~\ref{sec:algorithmic} we define uniform multicovers and use
these to define our procedural CAR model.  We show that it generates
all and only CAR mechanisms (Theorem~\ref{thm:algorithmic}). In
Section~\ref{sec:discussion} we discuss our CAR model in detail.  We
show (Theorem~\ref{thm:fibo}) that it gives rise to an exponential
lower bound on the height of the multicovers needed in
Theorem~\ref{thm:algorithmic}. The proofs are given in the final
section.

\section{Preliminaries}
\label{sec:prelim}

Let $E$ be a finite non-empty set, containing $n$ elements.  A
coarsening mechanism is a probabilistic rule which replaces any point
$x$ in $E$ with a subset $A$ of $E$ containing $x$.  Thus  a coarsening
mechanism is specified by a collection of (conditional) probabilities
$\pi_A^x$ such that for all $x$, $\sum_{A\ni x} \pi_A^x=1$.
Intuitively, $x$ is generated by some process which for simplicity we
will refer to as `Nature'. But rather than observing $x$
directly, the statistician observes a coarsening of $x$, i.e. a set
$A$ containing $x$. We call $x$ the {\em underlying outcome\/} and $A$
the corresponding {\em observation}. The coarsening mechanism
determines the $A$ that is observed given $x$; $\pi_A^x$ is the
probability of observing the set $A$ with $A \ni x$, given that Nature
has generated $x$.  We define the {\em support\/} of such a coarsening
mechanism as the set of $A \subseteq E$ for which
$\pi_A^x > 0$ for some $x \in E$.

A coarsening mechanism
satisfies the CAR (coarsening at random) property if and only if for
all $x,x'\in A$,
\begin{equation}
\label{eq:cardef}
\pi_A^x=\pi_A^{x'}=\pi_A \text{\ , say.} 
\end{equation} 
Intuitively, this means that the probability of observing $A$ is the
same for all $x$ that are contained in $A$: the coarsening is done `at
random', independently of the underlying $x$.  We note that
(\ref{eq:cardef}) is the definition of CAR employed by
\cite{GillLR97}. It is called ``strong CAR'' by \cite{Jaeger05a}.
The definition is explained in detail by
\cite{GillLR97} and \cite{Jaeger05a}; motivation, practical
relevance and applications of the CAR property are discussed
extensively by \cite{GillLR97} and \cite{GrunwaldH03}.

Definition (\ref{eq:cardef}) shows that a CAR mechanism is specified by a
collection of probabilities $\pi_A$ indexed by the nonempty subsets
$A$ of $E$ satisfying
\begin{equation}
\sum_{A\ni x} \pi_A~=~1\quad\forall x\in E.
\end{equation}
We can therefore represent a CAR mechanism by the vector
$\boldpi=(\pi_A: \emptyset\subset A \subseteq E)$, where we assume
the subsets $A$ to be ordered in some standard manner.  For a given
finite set of CAR mechanisms $\boldpi_1, \ldots, \boldpi_p$, and
any probability vector $\boldlambda = (\lambda_1, \ldots,
\lambda_p)$, we define their {\em mixture\/} $\boldpi' = \lambda_1
\boldpi_1 + \ldots + \lambda_p \boldpi_p$.  The following two
observations are immediate:
\begin{enumerate}
\item For each partition of $E$, there is a unique CAR mechanism that has exactly that partition as its support (for each set $A$ in the partition, $\pi^x_A = \pi_A = 1$, for all $x \in A$). 
\item Each finite mixture of CAR mechanisms again represents a CAR mechanism.
\end{enumerate}
These two observations suggest a simple procedural CAR model: Fix some
integer $p > 0$ and pick $p$ (arbitrary) partitions ${\cal E}_1, \ldots, {\cal
  E}_p$ of $E$. Each of these induces a unique corresponding CAR
mechanism. Now fix an arbitrary distribution $\boldlambda =
\lambda_1, \ldots, \lambda_p$ on ${\cal E}_1, \ldots, {\cal E}_p$.
The coarsened data are now generated by first, independently of the
underlying $x$, selecting one of the $p$ partitions according to the
distribution $\boldlambda$. Then, within the chosen partition, the
unique $A$ is generated which contains the underlying $x$. One can think
of each partition as a `sensor' with the help of which the data are
observed. The procedure amounts to selecting a sensor completely at
random, independently of the underlying $x$ generated by Nature. This
procedural CAR model is called the {\sc CARgen} procedure by
\cite{GrunwaldH03}. The `parameters' of this procedure are the
number of partitions $p$, the
partitions ${\cal E}_1, \ldots, {\cal E}_p$ and the distribution
$\boldlambda$.  Clearly, for every setting of the parameters, the
resulting algorithm defines a CAR mechanism. One may be tempted to
think that, by an appropriate setting of the parameters, {\em all\/}
CAR mechanisms can be simulated by {\sc CARgen}, but the following
example shows that this is not the case:
\begin{example} \rm  {\bf \citep{GillLR97}}
\label{ex:boem}
Let $E = \{1,2,3\}$, $A_{12} = \{1,2 \}, A_{23} = \{2, 3\}$ and $A_{31} = \{3, 1 \}$. Consider the coarsening mechanism $\boldpi^*$ defined by 
\begin{equation}
\label{eq:boem}
\pi^{*1}_{A_{12}} = \pi^{*2}_{A_{12}} =
\pi^{*2}_{A_{23}} = \pi^{*3}_{A_{23}} = \pi^{*3}_{A_{31}} = \pi^{*1}_{A_{31}} =
\frac{1}{2},
\end{equation}
and $\pi^{*x}_A = 0$ for all other $x \in E, A \subseteq E$.  By
(\ref{eq:cardef}) it is immediately seen that this is a CAR mechanism.
But because the support of the mechanism is not a union of partitions
of $E$, it cannot be simulated by the {\sc CARgen} procedure.
\end{example}
The example shows that the {\sc CARgen} procedure is incomplete: there
exist CAR mechanisms which cannot be represented by any parameter
setting of {\sc CARgen}. The
question is now whether there exist `natural' procedural CAR models
which are complete. In previous work, two candidates for such models
were proposed: Gr\"unwald and Halpern's (\citeyear{GrunwaldH03}) {\sc
  CARgen}$^*$ (an extension of {\sc CARgen} described above) and
Jaeger's (\citeyear{Jaeger05b}) {\em Propose-and-Test\/}-model. Both of
these suffer from the frailty property mentioned in the introduction:
rather than producing CAR mechanisms for all parameter settings, the
parameters need to be fine-tuned.  In previous work, one other
procedural model has been proposed which, like {\sc CARgen}, produces
CAR mechanisms for all settings of its parameters. However, as shown
by \cite{Jaeger05b}, this {\em randomized monotone coarsening\/}
model \citep{GillLR97} is in fact equivalent to {\sc CARgen}: both
can simulate exactly the set of `coarsening completely at random'
(CCAR) mechanisms.  In fact, \cite[Theorem 4.17]{Jaeger05b} shows that
any CAR mechanism that is not CCAR is, in a certain sense, nonrobust.
For the details of Jaeger's definition of robustness we refer to
\citep{Jaeger05b}. Briefly, he supposes that a CAR mechanism
involves an auxiliary randomization, and defines robustness in
terms of robustness to changes in the distribution of the
auxiliary variable.

Jaeger's result suggests that there exists no procedural CAR model
that is both complete and does not require any parameter tuning. Yet
in Section~\ref{sec:algorithmic}, we exhibit a simple extension of the
{\sc CARgen} procedure which achieves exactly this, as long as we are
able to sample from a uniform distribution. The procedural model will
be based on a geometric interpretation of CAR which we present below.
\section{A Geometric View of CAR}
\label{sec:geometric}
We have already indicated that a finite mixture of CAR mechanisms
$\boldpi$ is itself a CAR mechanism. Hence, for a given finite sample
space $E$ the set of all CAR mechanisms defined with respect to $E$
forms a convex body in Euclidean space. In Theorem~\ref{thm:geometric}
we show that this body is a polytope with a finite number of extreme
points, the vertices of the polytope. 
In order to characterize these extreme points, we first note that the
support of a CAR mechanism is always a cover of
$E$. With any cover of $E$ we associate its incidence matrix: the
matrix $M$ with rows indexed by $x\in E$, columns indexed by $A$ in
the support, and elements $\mathbb 1_{\{x\in A\}}$. An incidence
matrix of a cover is a matrix of $0$'s and $1$'s with at least one $1$
in every row and column.
We now use these incidence matrices to define extreme CAR mechanisms
in an algebraic way. Theorem \ref{thm:geometric} below states that
these CAR mechanisms are also extreme points in the geometric sense,
justifying our terminology.

In the sequel,  vectors are always column vectors, 
even if we lazily list the elements in a row.
$\vec{0}$ and $\vec{1}$ denote
vectors of $0$'s and $1$'s respectively, whose length
depends on the context.

Take the incidence matrix $M$ of an arbitrary cover $(A_1, \ldots,
A_m)$ of $E$.  If the equation $M\vec{z}=\vec{1}$ has a nonnegative
solution, then this solution $\vec{z} = (z_1, \ldots, z_m)$ represents
a CAR mechanism $\boldpi$, where for any $A_j$ appearing in the cover,
$z_j = \pi_{A_j}$, and for any $A$ not appearing in the cover, 
$\pi_A = 0$ (see also \cite{GrunwaldH03}, who explain this in detail).
We call $\boldpi$ a CAR  mechanism corresponding to $M$. 
\begin{definition}\label{def:extreme}
\rm 
We call $\boldpi$ an {\em extreme CAR mechanism\/} if it corresponds to
an incidence matrix $M$ of a cover $(A_1, \ldots, A_m)$ such
that $M \vec{z} = \vec{1}$ has a unique, and strictly positive,
solution. 
\end{definition}
By definition, a CAR mechanism is extreme if and only if it is the 
\emph{only} CAR mechanism  with the same support.
It is easily checked that the mechanism $\boldpi^*$ of
Example~\ref{ex:boem} is an example of an extreme CAR mechanism: it is
the only CAR mechanism with support $A_{12}, A_{23}, A_{31}$.
The uniqueness also implies
that the support of an extreme CAR mechanism cannot have more than $n$
elements (the size of $E$). It is clear that the 
number of extreme CAR mechanisms, for
given $E$, is finite. We can find them all by enumerating and testing
all covers of $E$ with $m \le n$ elements.
\begin{theorem}
\label{thm:geometric}
Every CAR mechanism is a mixture of extreme CAR mechanisms.
\end{theorem}
In other words,  all CAR mechanisms can be
represented by randomly choosing, independently of $x$, one of a finite
set of 
extreme CAR mechanisms. In the next section, we show that all such
extreme mechanisms are of a simple and natural form. This will lead to
Theorem~\ref{thm:algorithmic}, a direct corollary of
Theorem~\ref{thm:geometric}, giving an algorithmic characterization of
CAR. 
\section{An Algorithmic View of CAR}
\label{sec:algorithmic}
Our procedure is based on the notion of a {\em uniform multicover},
which we now define. A \emph{$k$-multicover} of $E$, or just $k$-cover
for short, is a collection of nonempty subsets of $E$, allowing
multiplicities, such that for each $x\in E$, precisely $k$ of the sets
(some of which may be the same) contain $x$.  Thus a $1$-cover is an
ordinary partition of $E$.  By a \emph{uniform multicover} we mean a
$k$-cover for some $k\ge 1$. The \emph{height} of a uniform multicover
is its value of $k$.  The {\em support\/} of a multicover is the set
of subsets of $E$ in the multicover.

A $k$-cover is specified by its support and by the multiplicity of
each set in its support. Thus, to each nonempty subset $A$ of $E$ there
corresponds a nonnegative integer $n_A$ such that $n_A=0$ if $A$ is absent from
the $k$-cover,  otherwise $n_A>0$ is the multiplicity of $A$ in the
$k$-cover. The $n_A$ have to satisfy
\begin{equation}
\label{eq:treini}
\sum_{A\ni x} n_A~=~k\quad\forall x\in E.
\end{equation}
For a given $k$-cover we can now define a CAR mechanism by setting
\begin{equation}
\label{eq:kbound}
\pi_A~=~n_A/k\quad\forall A\subseteq E.
\end{equation}
The algorithmic interpretation is as follows: Nature generates some $x
\in E$. The coarsening mechanism investigates which $A$ in the uniform
multicover contain $x$. There are exactly $k$ such $A$, including
multiplicities, whatever $x$. We choose one of these uniformly at
random, i.e. each $A$ with $x \in A$ is chosen with probability $1/k$.

Conversely, any CAR mechanism for which all the CAR probabilities
$\pi_A$ are rational numbers is generated by a 
$k$-cover with $k$ equal to the lowest common
multiple of the denominators of the $\pi_A$. 
We call CAR mechanisms obtained in this way \emph{rational}.
The \emph{rational CAR mechanisms are precisely 
the CAR mechanisms generated by a uniform multicover}.
Note that if $k$ and all $n_A$ share a common factor, we can
divide by this factor without changing the $\pi_A$. We consider
such multicovers as equivalent and take the multicover
with the smallest $k$ as representative of the class. In this way, each
rational CAR mechanism corresponds to exactly one uniform multicover, and vice
versa. 
We can make the connection to Theorem~\ref{thm:geometric} by noting that
\begin{fact}
  Every extreme CAR mechanism is rational. Thus, it is generated by a
  uniform multicover.
\end{fact}
This follows directly from the fact that the matrix $M$ in
Definition~\ref{def:extreme} is a $0/1$-matrix and the solution of $M
\vec{z} = \vec{1}$ is unique.

As stated above, for each rational CAR
mechanism there is a unique uniform multicover which generates it. We
can thus define an ``extreme multicover'' as a
uniform multicover that generates an extreme CAR mechanism. Using
Theorem 1, it is
easily shown that extreme multicovers are just those uniform multicovers 
that do not contain a subset that is also a uniform multicover (we
omit the details of the reasoning). 

We may now define a procedural CAR model by first fixing a finite
number $p$ of arbitrary uniform multicovers ${\cal C}_1, \ldots, {\cal
  C}_p$. We then fix an arbitrary distribution $\boldlambda =
(\lambda_1, \ldots, \lambda_p)$ on ${\cal C}_1, \ldots, {\cal C}_p$.
The coarsened data are now generated by first, independently of the
underlying $x$, selecting one of the $p$ uniform multicovers according
to the distribution $\boldlambda$. Suppose we have chosen multicover
${\cal C}_j$ with height $k_j$.  Then among the $k_j$ sets in ${\cal
  C}_j$ which contain $x$, we choose one uniformly at random, with
probability $1/k_j$.  This procedural CAR model is a simple extension
of {\sc CARgen} (Section~\ref{sec:prelim}), where the role of
partitions is taken over by the more general uniform multicovers.
Like {\sc CARgen}, it simulates a CAR mechanism for all parameter
settings; no fine-tuning is needed.  Theorem~\ref{thm:algorithmic}
below, part 2 (a corollary of Theorem~\ref{thm:geometric}) states that
by appropriately setting the parameters, we can simulate {\em all\/}
CAR mechanisms. Before presenting the theorem, we continue our
example.
\begin{example} \rm {[\bf Example~\ref{ex:boem} continued]}
  The collection ${\cal C} = \{A_{12}, A_{23}, A_{31} \}$ is a uniform
  multicover of $E$ with height 2. Consider a simple instantiation of
  the procedural CAR model we described above, with just one
  multicover ${\cal C} = {\cal C}_1$, so that $\boldlambda = (1)$.
  For each $x$ chosen by Nature, there will be exactly two elements of
  ${\cal C}$ which contain $x$. We select between these with
  probability $1/2$. It is immediately clear that this algorithm
  simulates the CAR mechanism $\boldpi^*$ described in
  Example~\ref{ex:boem}. An implementation of this
  mechanism requires a fair coin toss. If the coin is biased the
  CAR property can be lost. 
Relatedly,
the mechanism is not robust
  in Jaeger's sense.
\end{example}

\begin{theorem}
\label{thm:algorithmic} \ \\
\begin{itemize}
\item[1.]{Every CAR mechanism can be arbitrarily well approximated by a
    rational CAR mechanism, i.e. for all CAR mechanisms $\boldpi$, all $\epsilon >
    0$, there exists a rational CAR mechanism $\boldpi'$ such that $\| \boldpi -
   \boldpi' \| < \epsilon$.}
\item[2.]{Every CAR mechanism is exactly equal to a finite mixture of
    extreme (and hence rational) CAR mechanisms.}
\end{itemize}
\end{theorem}
We extensively discuss this theorem in the next section.

\section{Discussion}
\label{sec:discussion}
Theorem~\ref{thm:algorithmic} shows that there is an easy
probabilistic algorithm which approximates each CAR mechanism
arbitrarily well, and that a randomized version of the algorithm
reproduces each one exactly. Since the rational numbers form a dense
subset of the reals, part 1 of Theorem~\ref{thm:algorithmic} is, in a
sense, trivial. The real innovation is part 2, which shows that each
CAR distribution can be represented {\em exactly\/} as a mixture of a
finite set of candidate rational mechanisms.

 No fine tuning of parameters is required
to ensure the CAR properties so the algorithms do have a robustness
property. We just need to be able to choose uniformly at random from a
finite set.  Of course, if one perturbs the uniform distribution over
the $k$ sets containing a point $x$, one will in general destroy the
CAR property -- this is the reason that our result does not contradict
Jaeger's (\citeyear{Jaeger05b}) Theorem 4.17. 
For this reason, some readers may not want to call the procedure
`robust'.  However, the (weaker) claim that the algorithm requires no
parameter tuning seems indisputable: we can hardly think of
implementing a uniform distribution as `parameter tuning'. Unlike the
parameters in earlier complete procedural CAR models, which could vary
from situation to situation and were hard to determine, the uniform
distribution is universal and easy to determine. If the device we use
to generate a uniform distribution does not work perfectly, our
procedural model will slightly violate CAR, hence one might perhaps
say it is `nonrobust'; but devices used to generate a uniform
distribution (coins, dices) exist, and usually do not arise as
fine-tuned versions of devices that can generate a whole range of
distributions; hence one cannot say that our model requires `fine
tuning'.

The reason that earlier complete procedural CAR models did require
parameter tuning, was that their parameters had to satisfy complicated
constraints (see, for example, Example 4.7 in \citep{Jaeger05b}). As
remarked by M. Jaeger, we do pay a price for avoiding these parameter
constraints: we now have complicated constraints (\ref{eq:treini}) on
multiplicities of sets appearing in multicovers. Such constraints are
arguably more natural than constraints on continuous-valued
parameters, at least as long as the multicovers involved are not too
complex.  Unfortunately, in order to span all CAR mechanisms, we
sometimes need highly complex multicovers, as we show below. 
This limits the importance of our procedural model, as we discuss further below.

We can measure the complexity of multicovers in terms of their height.
Since the row rank of $M$ equals its (full) column rank, $m$, we can
delete rows obtaining an $m\times m$ nonsingular matrix $M_0$.
Deleting the corresponding rows from $\vec 1$ also, we obtain $\vec{z}
= M_0^{-1} \vec{1}$. It follows by the standard expression of matrix
inverse in terms of determinants that the value of $k$ appearing in
(\ref{eq:kbound}) is bounded by $m!$.
Hence, the height of the extreme multicovers that can be defined on a sample
space of size $|E| = n$ is upper bounded by $n!$. But is this too
pessimistic? 
Unfortunately not, or at least, not significantly:
our next and last
theorem gives an exponential \emph{lower} bound on the \emph{maximal} height of an
extreme multicover. It turns out that
this grows at least as fast as the celebrated Fibonacci
numbers, defined as  $F_1 = 1, F_2 = 1,$ and for $j \geq 3$, $F_j = F_{j-1} +
F_{j-2}$.

Theorem~\ref{thm:fibo} below considers
$n \times n$ matrices $S_n$ inductively defined as follows: $S_1 =
(1)$. For odd $n$, $S_{n+1}$ is constructed from $S_n$ by setting
$$
S_{n+1} = \left( \begin{array}{cc}
1 & \vec{0}^\top \\
\vec{0} & S_n 
\end{array}
\right).
$$
For even $n$, $S_{n+1}$ is constructed from $S_n$ by setting
$$
S_{n+1} = \left( \begin{array}{cc}
0 & \vec{1}^\top \\
\vec{1} & S_n 
\end{array}
\right).
$$
This is easier than it seems: the pattern should be obvious from the
example $n = 9$, shown in Figure~\ref{fig:lex}.
\begin{theorem}
\label{thm:fibo}
For odd $n > 0$, the equation $S_n \vec{z} = \vec{1}$ has the unique
solution
$$\vec{z} = \left( \frac{F_{n-1}}{F_n}, \frac{F_{n-2}}{F_n} \ldots,
  \frac{F_2}{F_n}, \frac{F_1}{F_n},\frac{1}{F_n} \right),$$
so that
$S_n$ represents an extreme point for sample spaces with size $|E| =
n$, with height $k = F_n$.
\end{theorem}
\begin{figure}[t]
{\small  
$$ \left( {\tt \begin{array}{ccccccccc}
0 & 1 & 1 & 1 & 1 & 1 & 1 & 1 & 1 \\
1 & 1 & 0 & 0 & 0 & 0 & 0 & 0 & 0 \\
1 & 0 & 0 & 1 & 1 & 1 & 1 & 1 & 1 \\
1 & 0 & 1 & 1 & 0 & 0 & 0 & 0 & 0 \\
1 & 0 & 1 & 0 & 0 & 1 & 1 & 1 & 1 \\
1 & 0 & 1 & 0 & 1 & 1 & 0 & 0 & 0 \\
1 & 0 & 1 & 0 & 1 & 0 & 0 & 1 & 1 \\
1 & 0 & 1 & 0 & 1 & 0 & 1 & 1 & 0 \\
1 & 0 & 1 & 0 & 1 & 0 & 1 & 0 & 1  
\end{array} }
\right)$$}
\caption{\label{fig:lex} The matrix $S_9$, an example of the
  matrices $S_n$ figuring in Theorem~\ref{thm:fibo}.
}
\end{figure}
The theorem implies that the maximal height of an extreme
multicover grows exponentially fast
with $n$; also, the maximal needed multiplicity of a set in 
an extreme multicover  grows exponentially fast
with $n$. We interpret this result as follows. 

Uniform multicovers are important in two ways:
\begin{enumerate}
\item They lead to an attractive algorithmic
characterization of CAR that requires no fine-tuning of parameters
(Theorem~\ref{thm:algorithmic}). 
\item They induce a hierarchy of CAR models that could be of use in
  statistical applications. We elaborate on this below.
\end{enumerate}
Yet apart from these applications,
the importance of uniform multicovers in understanding CAR is limited -- the maximal needed
height of the multicover grows exponentially fast with $n$, so though
the idea of the algorithm is simple, its detailed specification can be complex.
Thus, we can neither say that our characterization provides a truly simple
description of every CAR mechanism, nor that our multicover CAR mechanisms
always correspond to some `natural' process. While it seems reasonable to
suppose that low-height multicovers may be good models for some
processes occurring in nature, the same cannot be said for
exponentially high multicovers, and our Theorem~\ref{thm:fibo} does show
that we need to take these into account.

Jaeger's (\citeyear{Jaeger05b}) robustness Theorem 4.17 suggests that
the CAR mechanisms occurring in nature are those generated by
randomized 1-covers. Our characterization nuances this somewhat,
suggesting that in some situations $k$-covers for small $k > 1$ may
also be reasonable models. Indeed, the hierarchy of CAR mechanisms
induced by our algorithm suggests a statistical estimation procedure
for parsimoniously estimating CAR mechanisms and their parameters.
Such a procedure would penalize the fit of a proposed CAR mechanism to
the data. The penalization would be some function of the number of
extreme multicovers needed to express the mechanism, and the height of
each of these.  Alternatively one could use just one multicover, not
necessarily extreme, and penalize its height.  This could be done
either explicitly, by adding a regularization term to the likelihood,
or implicitly, by the use of suitable Bayesian priors.

Such procedures could be useful in practice if one seriously believed
that the data is CAR but quite possibly, not CCAR.  One could hope in
this way to combine the advantages of asymptotic validity and even go
for asymptotic efficiency, with good small sample behavior. However,
our results can also be read in a different way. Though we found an
appealing way to model CAR, it remains the fact that there do not seem
to be so many good reasons in practice, in general, to assume CAR but
not CCAR.  Therefore, if one is prepared to assume CAR, one is likely
to be also prepared to assume CCAR. Though the distinction concerns a
``nuisance'' part of the model, and indeed, in likelihood approaches
is invisible by the likelihood factorization implied by CAR, one can
capitalize on the extra knowledge for instance in order to obtain
better small sample properties of estimators, at the cost of loss of
asymptotic efficiency.

A final view is that the extra generality obtained by relaxing CCAR to
CAR is illusory. If one does not believe in CAR, one has no option but
to start modelling and estimating the coarsening mechanism. Jaeger
(\citeyear{Jaeger06b,Jaeger06a}) has made some proposals in this
direction which seem promising. Another possibility, so far not
explored, is to use the notion of \emph{relative} rather than
\emph{absolute} CAR introduced by \cite{GillLR97}.  The point of CAR
is that, in likelihood inference, one can analyse coarsened data
\emph{as if} the coarsening mechanism had been fixed in advance as any
particular CAR mechanism, and specifically therefore, as if coarsening
by an independently fixed-in-advance partitioning of the sample space.
Relative CAR means CAR relative to some other specific (non CAR)
coarsening mechanism: the likelihood factors; the interesting part is
the same \emph{as if} the data had been coarsened by the reference
coarsening model; the nuisance part can be used for inference
concerning \emph{which} coarsening mechanism has generated the data,
out of the mechanisms in the family implied by the reference
mechanism. It would be interesting to explore this possibility in more
detail. 

\section{Proofs}
\subsection{Proof of Theorem~\ref{thm:geometric}}
We show below  
that the set of all CAR mechanisms forms a convex polytope and
characterize the extreme points in terms
of linear algebra, corresponding to Definition \ref{def:extreme}.

A CAR mechanism is a collection of numbers $\pi_A$ indexed by the
nonempty subsets $A$ of a finite set $E$. They must satisfy two sets
of constraints: the inequalities $\pi_A\ge 0$ for each $A$, and the
equalities $\sum_{A\ni x}\pi_A=1$ for each $x$, both of which are
obviously linear. Together the constraints imply that $\pi_A\le 1$ for
all $A$.  Collecting the $\pi_A$ into a vector $\boldpi$ we see that the
set of all $\boldpi$ is a convex, compact polytope since it is bounded and
is the intersection of a finite number of closed half-spaces
(one for each inequality constraint) and hyperplanes (one for each
equality constraint).
Hence each $\boldpi$ is a convex combination of the extreme points of the
polytope, of which there are a finite number in total.

The polytope lives in the affine subspace of all vectors 
$\boldpi$ satisfying the equality
constraints $\sum_{A\ni x}\pi_A=1$ for each $x$. 
Since $\boldpi$ has $2^{|E|}-1$ components (the number of nonempty
subsets of $E$) and there are $|E|$ constraints, it follows that  the
dimension of this affine subspace is $2^{|E|}-1 - |E|$. The polytope is
just the intersection of that affine subspace with the positive orthant.
Within the 
affine subspace, each face of the polytope corresponds to
one of the hyperplanes $\pi_A= 0$. Each vertex of the polytope
is the unique meeting point of a number of faces; one for
each $A$ such that $\pi_A=0$.
Thus to each vertex is associated a collection of subsets $A$
such that if we set the corresponding $\pi_A$ equal to $0$ in the
equations $\sum_{A\ni x}\pi_A=1$ for all $x$, there is a unique and
strictly positive solution in the remaining $\pi_A$.
Conversely, any such collection of $A$ defines a vertex.

The subsets $A$ \emph{not} in the collection define the support
of the 
extreme
CAR mechanism $\boldpi$ under consideration. 
Let $M$ be its incidence matrix:  the
matrix of zeros and ones with rows indexed by elements $x\in E$, columns indexed by 
$A$ in the support, and with entries $\mathbb 1_{\{x\in A\}}$.
Write $\boldpi_{0}$ for the vector of $\pi_A$ for $A$ in the support. 
In matrix form, the equations which must have a unique and positive
solution $\vec z=\boldpi_0$ can be written
\begin{equation}
M\vec z = \vec{1},
\end{equation}
and we have proved that there is a one-to-one correspondence between
vertices of the polytope and incidence matrices $M$ of covers of $E$
such that this equation has a unique and positive solution.
As we argued in Section~\ref{sec:algorithmic}, if the solution is unique it
has to be rational.

Combining these facts, extreme points of
the polytope of CAR mechanisms
correspond to covers of $E$ whose incidence
matrix $M$ is such that  $M \vec{z}=\vec{1}$ has
a unique solution, and the solution is
strictly positive. 
\paragraph{Remark} {A 
condition equivalent to $M \vec{z} = \vec{1}$ having a unique positive
solution (Farkas' lemma in the theory of linear
programming, \cite[Chapter 7]{Schrijver86}),
is that $M$ has full column rank, and, if $\vec y$ is such that 
(a) $\vec y^\top M\ge \vec{0}$,
then (b) $\vec y^\top\vec{1} \ge 0$, with equality in (b)
implying equality in (a). 
By arguments from integer programming (see
again \citep{Schrijver86}) one may restrict 
here to vectors $\vec y$ of integers. \cite{Jaeger05b} gives
a version of this condition for the existence of
a CAR mechanism with given support -- he does not
demand full rank since he does not ask for
uniqueness. Though more combinatorial
in nature, this version of the condition for
extremality does not seem to be much more
useful, except perhaps for helping one to
show that certain covers do \emph{not}
lead to solutions.}
\subsection{Proof of Theorem~\ref{thm:algorithmic}}
Theorem~\ref{thm:algorithmic} is, in fact, a direct corollary of Theorem~\ref{thm:geometric}. Namely, each extreme point is rational and therefore corresponds to a uniform multicover.
Every point in a polytope can be written as a mixture of its extreme points.
This gives us item 2.  Item 1 follows by
considering the rational convex combinations of the extremes, which
lie dense in all convex combinations.
\subsection{Proof of Theorem~\ref{thm:fibo}}
We prove the theorem by induction on $n$. For $n = 1$, the result
trivially holds. Now suppose the result holds for $S_{n-1}$, for some
even $n >1$. Thus,
$
S_{n-1} \vec{q} = \vec{1}
$
has a unique solution
\begin{equation}
\label{eq:start}
\vec{q} = (q_1, \ldots, q_{n-1}) = 
\left( \frac{F_{n-2}}{F_{n-1}}, \frac{F_{n-3}}{F_{n-1}} \ldots,
  \frac{F_2}{F_{n-1}}, \frac{F_1}{F_{n-1}},\frac{1}{F_{n-1}}
\right).
\end{equation} 
We prove the theorem by showing that this implies that
\begin{equation}
\label{eq:ind}
S_{n+1} \vec{r} = \vec{1} 
\end{equation}
has the unique solution:
\begin{equation}
\label{eq:finish}
\vec{r} = (r_1, \ldots, r_{n+1}) = 
\left( \frac{F_{n}}{F_{n+1}}, \frac{F_{n-1}}{F_{n+1}} \ldots,
  \frac{F_2}{F_{n+1}}, \frac{F_1}{F_{n+1}},\frac{1}{F_{n+1}}
\right).
\end{equation} 
To prove (\ref{eq:finish}),
note first that to each
row of
(\ref{eq:ind}) corresponds a linear equation. Writing the equations
corresponding to the first two rows explicitly and the equations
corresponding to rows 3 to $n+1$ in matrix form, and reordering terms, we see that
(\ref{eq:ind}) is equivalent to:
\begin{eqnarray}
\label{eq:r2rest} 
r_2 & = &  1- \sum_{i=3}^{n+1} r_i \\
\label{eq:r21}
r_2 & = & 1- r_1 \\
\label{eq:r1rest}
S_{n-1} (r_3, \ldots, r_{n+1})^T & = &  1 - r_1,
\end{eqnarray}
where, by our inductive assumption, the last equality implies 
\begin{equation}
\label{eq:r1all}
(r_3, \ldots, r_{n+1}) = (1 - r_1) (q_1, \ldots, q_i).
\end{equation}
and in particular
\begin{equation}
\label{eq:r1sum}
\sum_{i=3}^{n+1} r_i = (1 - r_1) \sum_{i=1}^{n-1} q_i .
\end{equation}
Combining (\ref{eq:r21}) with (\ref{eq:r2rest}), we get $r_1 =
\sum_{i=3}^{n-1} r_i$. Plugging this into (\ref{eq:r1sum}) gives
\begin{equation}
\label{eq:hallo}
\frac{r_1}{1-r_1} = \sum_{i=1}^{n-1} q_i
\end{equation}
where $q_i$ are given by (\ref{eq:start}).
We claim this has the unique solution $r_1 = F_n/F_{n+1}$. To see
this, note the following basic fact which follows immediately from
repeatedly substituting  the
definition $F_n = F_{n-1} + F_{n-2}$ on the left in (\ref{eq:fact}): 
\begin{fact}
For odd $n > 0$,
\begin{equation}
\label{eq:fact}
F_n = \sum_{i=1}^{n-2} F_j + 1.
\end{equation}
\end{fact}
The fact implies that the right-hand side of (\ref{eq:hallo}) is equal
to $F_n / F_{n-2}$. Plugging in our proposed solution $r_1 =
F_n/F_{n+1}$, the left-hand side of (\ref{eq:hallo}) becomes $F_n/(F_{n+1} - F_n) =
F_n/F_{n-2}$, so that (\ref{eq:hallo}) holds. This shows that
$r_1$ is indeed given by $F_n/F_{n+1}$. By (\ref{eq:r21}) it now
follows that $r_2 = F_{n-1}/F_{n+1}$, and, by (\ref{eq:r1all}), that
for $j \in \{3, \ldots, n+1\}$, $r_j = q_{j-2}/r_2 = s_j/F_{n+1}$,
where $(s_1, s_2, \ldots, s_{n-2}, s_{n-1}) = (F_{n-2}, F_{n-3}, \ldots, F_1,1)$.
This shows that (\ref{eq:finish}) is the unique solution of $S_{n+1} \vec{r}
= \vec{1}$, and thus completes the induction step.
The theorem is proved. 
\section{Acknowledgments}
We would like to thank Sasha Gnedin and Lex Schrijver for stimulating
conversations.  Lex Schrijver made some essential contributions to the
proof of Theorem~\ref{thm:fibo}. 

This work was supported by
the IST Programme of the European Community 
within the PASCAL Network of Excellence, 
IST-2002-506778. RDG's work was 
supported (during his previous affiliation) 
by the Dept.~of Mathematics, Utrecht University,
and, thanks to a visiting position, 
by the Thiele Centre, \AA rhus University. 
CWI is the Dutch national research 
institute for mathematics and computer
science. {\sc Eurandom} is funded by the Dutch science 
foundation, NWO, and Eindhoven University.

\bibliographystyle{chicago}

\bibliography{final}

\end{document}